\newcommand{\al}{\alpha}              \newcommand{\be}{\beta}
\newcommand{\ga}{\gamma}              \newcommand{\Ga}{\Gamma}
\newcommand{\de}{\delta}              
\newcommand{\lb}{\lambda}             \newcommand{\Lb}{\Lambda}
                  
\newcommand{\veps}{\varepsilon}       \newcommand{\vphi}{\varphi}

\newcommand{\calp}{{\mathcal P}}

\newcommand{\arch}{{\rm arch}}        \newcommand{\cl}{{\rm cl}}
        
      \newcommand{\Dom}{{\rm Dom}}

         \newcommand{\incl}{\subseteq}
             \newcommand{\sk}{\smallskip}
\newcommand{\es}{\emptyset}            \newcommand{\sm}{\setminus}
\newcommand{\limpl}{\Longrightarrow}
\newcommand{\lequi}{\Longleftrightarrow}
            \newcommand{\oo}{\infty}
           \newcommand{\wt}{\widetilde}

\newcommand{\barr}{\begin{array}}       \newcommand{\earr}{\end{array}}
\newcommand{\beq}{\begin{equation}}     \newcommand{\eeq}{\end{equation}}
\newcommand{\bit}{\begin{itemize}}      \newcommand{\eit}{\end{itemize}}
\newcommand{\blemma}{\begin{lemma}}     \newcommand{\elemma}{\end{lemma}}
\newcommand{\bproof}{\begin{proof}}     \newcommand{\eproof}{\end{proof}}
\newcommand{\bprop}{\begin{proposition}} \newcommand{\eprop}{\end{proposition}}
\newcommand{\btab}{\begin{tabular}}     \newcommand{\etab}{\end{tabular}}
\newcommand{\btheorem}{\begin{theorem}} \newcommand{\etheorem}{\end{theorem}}

\documentclass[reqno]{amsart}

\newtheorem{theorem}{\bf Theorem}

\newtheorem{lemma}{\bf Lemma}
\newtheorem{proposition}{\bf Proposition}

\begin{document}

\title
[NONLINEAR VERSIONS OF A VECTOR MAXIMAL PRINCIPLE]
{NONLINEAR VERSIONS \\ 
OF A VECTOR MAXIMAL PRINCIPLE}

\author{Mihai Turinici}
\address{
"A. Myller" Mathematical Seminar;
"A. I. Cuza" University;
700506 Ia\c{s}i, Romania
}
\email{mturi@uaic.ro}




\begin{abstract}
Some nonlinear extensions of the
vector maximality statement established by
Goepfert, Tammer and Z\u{a}linescu [Nonl. Anal., 39 (2000), 909-922]
are given. Basic instruments for these are the
Brezis-Browder ordering principle [Advances Math., 21 (1976), 355-364]
and a (pseudometric) version of it obtained in
Turinici [Demonstr. Math., 22 (1989), 213-228].
\end{abstract}

\maketitle

\section{Introduction}
\setcounter{equation}{0}

Let $Y$ be a (real) separated {\it locally convex space};
and $K$, some {\it (convex) cone} of it:
$\al K+\be K\incl K$, $\forall  \al, \be \in R_+:=[0,\oo[$.
In this case, the relation\  over $Y$\ 
\bit
\item[(a01)]
$(y_1,y_2\in Y$):\ \
$y_1\le_K y_2$\ if and only if\ \ $y_2-y_1\in K$
\eit
is reflexive and transitive; hence a {\it quasi-order};
in addition, it is {\it compatible}
with the linear structure of $Y$.
Let $H$ be another (convex) cone of $Y$ with $K\incl H$;
and pick some $k^0\in K\sm (-H)$.
Further, take some complete metric space $(X,d)$;
and introduce the quasi-order (on $X\times Y$)
\bit
\item[(a02)]
$(x_1,y_1) \succeq (x_2,y_2)$ iff $k^0 d(x_1,x_2)\le_K y_1-y_2$.
\eit
Finally, take some nonempty part $A$ of $X\times Y$.
For a number of both practical and theoretical reasons,
it would be useful to determine sufficient conditions
under which $(A,\succeq)$ has points with certain
maximal properties.
The basic 2000 result in the area obtained by
Goepfert, Tammer and Z\u{a}linescu \cite{goepfert-tammer-zalinescu-00},
deals with the
case $H=\cl(K)$ (=the {\it closure} of $K$).
Precisely, assume that
\bit
\item[(a03)]
$P_Y(A)$ is bounded below (modulo $K$):
$\exists \wt y\in Y$ with  $P_Y(A)\incl \wt y +K$
\item[(a04)]
if\ $((x_n,y_n))\incl A$ is ($\succeq$)-ascending
and\ $x_n \to x$\ then $x\in P_X(A)$  \\
and there exists\ $y\in A(x)$\ such
that $(x_n,y_n)\succeq (x,y)$, for all $n$.
\eit
[Here, for each $(x,y)\in A$, $A(x)$ (respectively, $A(y)$) stands for the
$x$-{\it section} (respectively, $y$-{\it section}) of (the relation) $A$;
and $P_X$, $P_Y$ are the {\it projection} operators from
$X \times Y$ to $X$ and $Y$ respectively].

\btheorem \label{t1}
Let the above conditions be in force.
Then, for each\ $(x_0,y_0)\in A$\ there exists\
$(\bar x,\bar y)\in A$ with
\beq \label{101}
\mbox{
$(x_0,y_0) \succeq (\bar x,\bar y)$ [hence $y_0\ge_K \bar y$]
}
\eeq
\beq \label{102}
\mbox{
$(\bar x,\bar y)\succeq (x',y')\in A$ imply $\bar x=x'$.
}
\eeq
\etheorem

This result includes the ones due to
Isac \cite{isac-96} and
Nemeth \cite{nemeth-86}
(which, in turn, extend
Ekeland's variational principle \cite{ekeland-79});
and the authors' argument is based on the Cantor intersection theorem.
Further, in his 2002 paper,
Turinici \cite{turinici-02}
proposed a different approach, via ordering principles related to
Brezis-Browder's \cite{brezis-browder-76}
(cf. Section 2); and stressed that, conclusions like before
are extendable to (non-topological) vector spaces  $Y$
under the choice
$H=\arch(K)$ (=the {\it Archimedean closure} of $K$).
It is our aim in this exposition to show that
a further enlargement of these facts is possible
(by the same techniques). This refers to the function
\bit
\item[(a05)]
$\Lb(t)=k^0 t,\ t\in R_+$ \quad (where $k^0$ is the above one)
\eit
being no longer linear; details will be given in
Section 4  (the Archimedean case) and
Section 5 (the non-Archimedean case).
The specific instrument of our investigations
(in this last circumstance) is the concept of {\it gauge}
function (developed in Section 3).
Finally, in Section 6, the relationships between our statement
and the recent variational principle in
Bao and Mordukhovich \cite{bao-mordukhovich-10}
are discussed.

\section{Brezis-Browder principles}
\setcounter{equation}{0}

{\bf (A)}
Let $M$ be some nonempty set. Take a
{\it quasi-order} $(\le)$
over $M$; as well as a function $\psi: M \to R_+$.
Call the point $z\in M$, $(\le,\psi)$-{\it maximal} when:
$w\in M$ and $z\le w$\ imply\ $\psi(z)=\psi(w)$.
A basic result about such points is the 1976
Brezis-Browder ordering principle \cite{brezis-browder-76}:

\bprop \label{p1}
Suppose that
\bit
\item[(b01)]
$(M,\le)$ is sequentially inductive:\\
each ascending sequence has an upper bound (modulo $(\le)$)
\item[(b02)]
$\psi$ is $(\le)$-decreasing ($x\le y \limpl \psi(x) \ge \psi(y)$).
\eit
Then, for each $u\in M$ there exists a $(\le, \psi)$-maximal
$v\in M$ with $u\le v$.
\eprop

This principle, including
Ekeland's \cite{ekeland-79},
found some basic applications to convex and nonconvex analysis
(cf. the above references).
So, a discussion about its key condition (b01)
would be not without profit.
Let $(Z,\le)$ be some quasi-ordered structure.
Take a function $\vphi: Z \to R\cup\{-\oo,\oo\}$;
and let $M$ be some nonempty part of $Z$.
For simplicity reasons, we let again $\vphi$ stand for
the restriction of $\vphi$ to $M$.
The following "relative" form of Proposition \ref{p1}
will be useful for us.

\bprop \label{p2}
Suppose (b02) holds (modulo $\vphi$), as well as
\bit
\item[(b03)]
$\vphi$ is inf-proper over $M$:\\ 
$\inf[\vphi(M)]>-\oo$\  and\  
$\Dom(\vphi):=\{x\in M; \vphi(x)< \oo\}\ne \es$
\item[(b04)]
$\Dom(\vphi)$ is sequentially inductive in $M$: each ascending\\
sequence in $\Dom(\vphi)$ is bounded above in $M$ (modulo $(\le)$).
\eit
Then, for each $u\in \Dom(\vphi)$ there exists $v\in \Dom(\vphi)$ with

{\bf i)} $u\le v$\ and\
{\bf ii)} $x\in M$, $v\le x$ imply $\vphi(v)=\vphi(x)$.
\eprop

\bproof
Let $u\in \Dom(\vphi)$ be arbitrary fixed.
Put $M(u,\le):=\{x\in M; u\le x\}$;
and introduce the function (from $M$ to $R_+$)
$\psi(x)=\vphi(x)-\vphi_*, x\in M$;
where $\vphi_*:=\inf[\vphi(M)]$.
By the imposed conditions, Proposition \ref{p1}
applies to $M(u,\le)$ and $(\le,\psi)$;
wherefrom the conclusion is clear.
\eproof

For the moment, Proposition \ref{p2} is a logical consequence of
Proposition \ref{p1}. The reciprocal is also true, by simply taking
$Z=M$, $\vphi=\psi$. Hence, these two results are logically
equivalent.
Note that the inf-properness condition (b03)
is not essential for the conclusion above
(cf. C\^{a}rj\u{a} and Ursescu \cite{carja-ursescu-93}).
Moreover, $(R,\ge)$ may be substituted by a separable ordering
structure $(P,\le)$ without altering the conclusion above;
see Turinici \cite{turinici-06} for details. 
Further aspects were discussed in
Altman \cite{altman-82}; see also
Kang and Park \cite{kang-park-90}.
\sk

{\bf (B)}
A semi-metric version of these developments may be given along the
following lines. Let $(M,\le)$ be taken as before.
By a {\it pseudometric} over $M$ we shall
mean any map $e:M\times M \to R_+$. If, in addition,
$e$ is
{\it reflexive} [$e(x,x)=0, \forall x\in M$],
{\it triangular} [$e(x,z)\le e(x,y)+e(y,z), \forall x,y,z\in M$] and
{\it symmetric} [$e(x,y)=e(y,x), \forall x,y\in M$],
we say that it is a {\it semimetric} (on $M$).
Suppose that we fixed such an object.
Call $z\in M$, $(\le,e)$-{\it maximal}, in case:
$w\in M$ and $z\le w$ imply $e(z,w)=0$.
[Note that, if (in addition) $e$ is
{\it sufficient} [$e(x,y)=0$ implies $x=y$],
this property  becomes;
$w\in M, z\le w \limpl  z=w$
(and reads: $z$ is strongly $(\le)$-maximal).
So, existence results of this type may be viewed as "metrical"
versions of the Zorn-Bourbaki principle].
To get such points, one may proceed as below.
Call the (ascending) sequence $(x_n)$ in $M$,  $e$-{\it Cauchy} when
[$\forall \de> 0, \exists n(\de)$: 
$n(\de)\le p\le q \limpl e(x_p,x_q)\le \de$];
and $e$-{\it asymptotic}, provided 
[$e(x_n,x_{n+1})\to 0$, as $n \to \oo$].
Clearly, each (ascending) $e$-Cauchy sequence is $e$-asymptotic too.
The reverse implication is also true when all such sequences are involved;
i.e., the global conditions below are equivalent:
\bit
\item[(b05)]
each ascending sequence is $e$-Cauchy
\item[(b06)]
each ascending sequence is $e$-asymptotic.
\eit
By definition, either of these will be referred to as $(M,\le)$
is {\it regular} (modulo $e$).
The following maximality result in
Turinici \cite{turinici-89}
is available.

\bprop \label{p3}
Assume that $(M,\le)$ is sequentially inductive and regular (modulo $e$).
Then, for each $u\in M$ there exists an $(\le,e)$-maximal
$v\in M$ with $u\le v$.
\eprop

This result includes the
Brezis-Browder ordering principle \cite{brezis-browder-76}
(Proposition \ref{p1}); to which it reduces in case
$e(x,y)=|\psi(x)-\psi(y)|$ (where $\psi$ is the above one).
The reciprocal inclusion is also true; we refer to
the quoted paper for details.

\section{Conical gauge functions}
\setcounter{equation}{0}

Let $Y$ be a (real) {\it vector space}.
Take a convex cone $L$  of $Y$ (cf. Section 1); which, in addition,
is {\it non-degenerate} [$L\ne \{0\}$] and {\it proper} [$L\ne Y$].
Denote by $(\le_L)$ its induced quasi-order (cf. (a01));
when $L$ is understood, we indicate this as $(\le)$, for simplicity.
Further, let the map $\Lb:R_+\to L$ be {\it normal} (modulo $L$):
\bit
\item[(c01)]
$\Lb(0)=0$ and $\Lb$ is strictly increasing (modulo $L$):\\
$\Lb(\tau)-\Lb(t)\in L\sm (-L)$,\ \ whenever $\tau> t$
\item[(c02)]
$\Lb$ is sub-additive:
$\Lb(t_1+t_2)\le \Lb(t_1)+\Lb(t_2)$,\
$\forall t_1,t_2\in R_+$.
\eit
Note that, as a consequence of this,
\beq \label{301}
\mbox{
$\Lb$ is strictly positive (modulo $L$):
$\Lb(t)\in L\sm (-L), \forall t\in R_+^0:=]0,\oo[$
}
\eeq
\beq \label{302}
\mbox{
$\Lb$ is subtractive:
$\Lb(t_1-t_2)\ge \Lb(t_1)-\Lb(t_2)$,\
$\forall t_1,t_2\in R_+,\ t_1\ge t_2$.
}
\eeq
Having these precise, denote (for $y\in Y$)
\bit
\item[(c03)]
$\Ga(L;\Lb;y)=\{s\in R_+; \Lb(s)\le y\},\ \
\ga(L;\Lb;y)=\sup \Ga(L;\Lb;y)$.
\eit
(By convention, $\sup(\es)=-\oo$).
We therefore defined a couple of functions $\Ga(.):=\Ga(L;\Lb;.)$
and $\ga(.):=\ga(L;\Lb;.)$ from $Y$ to $\calp(R_+)$ and
$R\cup \{-\oo,\oo\}$ respectively; the latter of these will be
referred to as the {\it gauge} function attached to $(L;\Lb)$.
For the particular case of {\it linear} normal functions
[i.e., the one of (a05), with $k^0\in L\sm (-L)$],
such objects were introduced (in the same context) by
Turinici \cite{turinici-02};
and these, in turn, appear as non-topological extensions of the
locally convex ones in
Goepfert, Tammer and Z\u{a}linescu \cite{goepfert-tammer-zalinescu-00}.
The present developments may therefore be viewed as  "nonlinear"
extensions of the preceding ones.
\sk

{\bf i)}
To begin with, note that 
for each $y\in L$, $\Ga(y)$ is {\it hereditary} 
($s\in \Ga(y)\limpl [0,s]\incl \Ga(y)$).
In addition, we have (by definition)
$y\in L \lequi  \Ga(y)\ne \es \lequi \ga(y)\in [0,\oo]$) 
(or, equivalently: $y\notin L \lequi \Ga(y)=\es \lequi \ga(y)=-\oo$).

{\bf ii)} 
The gauge function is increasing
[$y_1,y_2\in Y$, $y_1\le y_2$ $\limpl$ $\ga(y_1)\le \ga(y_2)$].

{\bf iii)}
Further,  $\ga$ is super-additive and subtractive:
\beq \label{303}
\mbox{
$\ga(y_1+y_2)\ge \ga(y_1)+\ga(y_2)$,\ \
whenever the right member exists
}
\eeq
\beq \label{304}
\mbox{
$\ga(y_1-y_2)\le \ga(y_1)-\ga(y_2)$,\ \
if $\ga(y_2)$=finite (hence $0\le \ga(y_2)< \oo$).
}
\eeq
Clearly, it will suffice proving the former one.
Without loss, assume that
$\ga(y_1)> 0$, $\ga(y_2)> 0$. 
By definition (and the hereditary property of $\Ga$)
$y_1\ge \Lb(t_1),\ y_2\ge \Lb(t_2)$,
whenever $0\le t_1< \ga(y_1),\ 0\le t_2< \ga(y_2)$;
so, combining with (c02), yields (for all such $(t_1,t_2)$)
$y_1+y_2\ge \Lb(t_1)+\Lb(t_2) \ge \Lb(t_1+t_2)$;
that is, $\ga(y_1+y_2)\ge t_1+t_2$.
This, and the arbitrariness of the precise couple, ends the
argument.

{\bf iv)} 
Finally, the {\it identity} relation is available:
\beq \label{306}
\mbox{
$\ga(\Lb(t))=t$, \quad for each $t\in R_+$.
}
\eeq
In fact, let $t\in R_+$ be arbitrary fixed; it will suffice
verifying that $\Ga(\Lb(t))=[0,t]$. Suppose not:
there exists $\tau> t$ with $\tau\in \Ga(\Lb(t))$. By definition,
$\Lb(\tau)\le \Lb(t)$; wherefrom $\Lb(\tau)-\Lb(t)\in (-L)$;
in contradiction to (c01).
As a consequence, $\ga$ is {\it proper}; i.e., 
$\Dom(\ga):=\{y\in Y; \ga(y)< \oo\}$ is nonempty. 
Moreover, $\Dom_L(\ga):=\Dom(\ga)\cap L$
is nonempty too; and we have the decomposition 
$\Dom(\ga)=(Y\sm L)\cup \Dom_L(\ga)$, with (in addition)
$\ga(Y\sm L)=\{-\oo\}$, $\ga[\Dom_L(\ga)]=R_+$; 
as results from (\ref{306}).
On the other hand, by super-additivity, we have the
{\it sup-translation} property:
$\ga(y+\Lb(t))\ge \ga(y)+t$, $\forall y\in Y$, $\forall t\in R_+$.
This inequality may be strict; just take  $y=-\Lb(\tau), t=\tau$,
for some $\tau> 0$.
\sk

Concerning the effectiveness of such a construction,
call the function $\psi:R_+\to R_+$, {\it normal}, when
$\psi(0)=0$ and $\psi$ is strictly increasing (on $R_+$)
as well as sub-additive (see above).
Note that such functions exist; such as, e.g.;
$\psi(t)=t^\lb,\ t\in R_+$, for some $\lb\in ]0,1]$.
Suppose that $\{\psi_1,...,\psi_m\}$ are endowed with such
properties; and take some points $\{k^1,...,k^m\}$ in $L\sm (-L)$.
Then, the function
\bit
\item[(c04)]
($\Lb:R_+ \to L$):\
$\Lb(t)=k^1 \psi_1(t)+...+k^m \psi_m(t), \quad  t\in R_+$
\eit
is a normal one, in the sense of (c01)+(c02).
The obtained class of all these covers the
linear one (expressed via (a05)); when (as precise)
these developments reduce to the ones in
Turinici \cite{turinici-02}.
Further aspects involving the locally convex (modulo $Y$) case
(and the same linear setting) may be found in
Goepfert, Riahi, Tammer and Z\u{a}linescu
\cite[Ch 3, Sect 10]{goepfert-riahi-tammer-zalinescu-03}.
see also
Gerth (Tammer) and Weidner \cite{gerth(tammer)-weidner-90}.

\section{Main result (Archimedean case)}
\setcounter{equation}{0}

With these  preliminaries, we may now return to
the question of the introductory part.
Let $Y$ be a (real) vector space; and $K$, some (convex) cone
of it. Denote by $(\le_K)$ the induced quasi-order (cf. (a01)).
Let $H$ be another (convex) cone of $Y$ with
$K\incl H$; and the map $\Lb:R_+\to K$ be {\it almost normal}
(modulo $(K,H)$):
\bit
\item[(d01)]
$\Lb(0)=0$ and $\Lb$ is strictly increasing (modulo $(K,H)$):\\
$\Lb(\tau)-\Lb(t)\in K\sm (-H)$,\ \ whenever $\tau> t$
\item[(d02)]
$\Lb$ is sub-additive (modulo $K$):\\
$\Lb(t_1+t_2)\le_K \Lb(t_1)+\Lb(t_2)$,\
$\forall t_1,t_2\in R_+$.
\eit
Note that, as a consequence of this (cf. Section 3)
\beq \label{401}
\mbox{
$\Lb$ is strictly positive (modulo $(K,H)$):
$\Lb(t)\in K\sm (-H)$, $\forall t\in R_+^0$
}
\eeq
\beq \label{402}
\mbox{  \btab{l}
$\Lb$ is subtractive modulo $K$):\\
$\Lb(t_1-t_2)\ge_K \Lb(t_1)-\Lb(t_2)$,\
$\forall t_1,t_2\in R_+,\ t_1\ge t_2$.
\etab
}
\eeq
In addition, let $(X,d)$ be a metric space.
The relation over $X\times Y$
\bit
\item[(d03)]
$(x_1,y_1) \succeq (x_2,y_2)$ iff $\Lb(d(x_1,x_2))\le_K y_1-y_2$
\eit
is reflexive and transitive (by the properties of $\Lb$);
hence a quasi-order on it.
Finally, take some (nonempty) part $A$ of $X\times Y$.
As in Section 1, we are interested to determine sufficient
conditions under which $(A,\succeq)$ has points with
certain maximality properties.
Note that, in the linear case of (a05), this problem is just the one in
Turinici \cite{turinici-02};
which (under $d$=complete) has a positive
answer in the context of ((a04) and)
\bit
\item[(d04)]
$P_Y(A)$ is bounded below (modulo $H$):
$\exists \wt y\in Y$ with  $P_Y(A)\incl \wt y +H$.
\eit
So, it is natural asking whether similar conclusions are
retainable in our "nonlinear" setting.
Loosely speaking, these depend on
the ambient convex cone $H$ being or not {\it Archimedean}.
So, two alternatives are open before us.
\sk

In the following, we discuss the former of these,
based on $H$ being endowed with such a property
(cf. Cristescu \cite[Ch 5, Sect 1]{cristescu-77}):
\bit
\item[(d05)]
$h,v\in H$ and [$h\tau\le_H v, \forall \tau\in R_+^0$] imply $h\in H\cap (-H)$.
\eit
As we shall see, a positive answer is available under
\bit
\item[(d06)]
each $(\succeq)$-ascending $e$-Cauchy sequence
$((x_n,y_n))\incl A$\\
is bounded above in $A$ (modulo $(\succeq)$).
\eit
Here, $e$ stands for the semi-metric on $X\times Y$ introduced as
\bit
\item[(d07)]
$e((x_1,y_1),(x_2,y_2))=d(x_1,x_2),\ (x_1,y_1),(x_2,y_2)\in X\times Y$.
\eit
The  (first) main result of our exposition is

\btheorem \label{t2}
Let the assumptions (d04)-(d06) hold.
Then, for each starting $(x_0,y_0)\in A$ there exists
$(\bar x,\bar y)\in A$ with the properties (\ref{101}) and (\ref{102})
(written for our data).
\etheorem

The latter of the conclusions above reads (under the precise convention)
$$  \mbox{
$(\bar x,\bar y)\succeq (x',y')\in A$ $\limpl$
$e((\bar x,\bar y),(x',y'))=0$.
}
$$
This suggests us a possible deduction of Theorem \ref{t2}
from Proposition \ref{p3}.
To see the effectiveness of such an approach, we need an auxiliary fact.

\blemma \label{le1}
Let $((x_n,y_n))$ be an $(\succeq)$-ascending sequence in $A$:
\bit
\item[(d08)]
$\Lb(d(x_n,x_m))\le_K y_n-y_m$, \quad whenever $n\le m$.
\eit
Then, $(x_n)$ is $d$-Cauchy in $P_X(A)$; hence
$((x_n,y_n))$ is $e$-Cauchy in $A$.
\elemma

\bproof {\bf (Lemma \ref{le1})}
Suppose that this would be not valid;
i.e. (as $d$ is symmetric), there  must be some $\veps> 0$ in
such a way that
\bit
\item[(d09)]
for each $n$, there exists $m> n$ with $d(x_n,x_m)\ge \veps$.
\eit
Inductively, we may construct a subsequence $(u_n=x_{i(n)})$ of
$(x_n)$ with $d(u_n,u_{n+1})\ge \veps$, for all $n$.
This in turn yields, for the corresponding subsequence\ $(v_n=y_{i(n)})$
of $(y_n)$, an evaluation like:
$\Lb(\veps) \le_K \Lb(d(u_n,u_{n+1}))\le_K v_n-v_{n+1}$,
for each $n\ge 1$.
But then, in view of (d04), one derives:
$q \Lb(\veps) \le_H v_1-v_{q+1} \le_H v_1- \wt y$,
for each $q\ge 1$.
This, along with (d05), gives $\Lb(\veps)\in K\cap (-H)$; in
contradiction with (\ref{401}).
Hence, the working assumption (d09) cannot hold;
and the claim follows.
\eproof

\bproof {\bf (Theorem \ref{t2})}
Let\ $((x_n,y_n))$ be an
$(\succeq)$-ascending sequence in\ $A$.\
By  Lemma \ref{le1},
$((x_n,y_n))$ is an $e$-Cauchy sequence in $A$;
which tells us that $(A,\succeq)$ is regular (modulo $e$).
Moreover, by (d06),
$((x_n,y_n))$ is bounded above (modulo $(\succeq)$) in $A$;
wherefrom, $(A,\succeq)$ is sequentially inductive.
Summing up, Proposition \ref{p3} is applicable to
$(A,\succeq;e)$; so that (from its conclusion) each
$a_0=(x_0,y_0)$ in $A$ is majorized (modulo $(\succeq)$)
by some $(\succeq,e)$-maximal $\bar a=(\bar x,\bar y)$ in $A$.
This gives the conclusions (\ref{101})+(\ref{102}) we need;
and completes the argument.
\eproof

In particular, when $\Lb$ is taken as in (a05) (and $d$ is complete)
Theorem \ref{t2} is just the related statement in
Turinici \cite{turinici-02};
which, as precise there, incorporates the (locally convex) one in
Goepfert, Tammer and Z\u{a}linescu \cite{goepfert-tammer-zalinescu-00}
(Theorem \ref{t1}).
This inclusion seems to be strict; because the choice  (c04)
of $\Lb$ cannot be (completely) reduced to the linear one
(appearing in all these papers). Some related aspects may be found in
Tammer \cite{tammer-92}.

\section{A completion (non-Archimedean case)}
\setcounter{equation}{0}

Now, the key regularity assumption used in the result above
is (d05). So, it is natural to discuss the alternative of
this being avoided.
As we shall see below, a positive answer is still available; but we
must restrict the initial set $A$ in a way imposed by the
associated (to $H$) gauge function.
\sk

Let $Y$ be a (real) vector space; and $K$, some (convex) cone
in it. Denote by $(\le_K)$ its associated quasi-order;
and let $H$ be another (convex) cone of $Y$ with $K\incl H$.
We also take a map $\Lb:R_+\to K$; which  is supposed to be
almost normal (modulo $(K,H)$) in the sense of (d01)+(d02).
Clearly, it is also normal (modulo $H$); so, we may construct
the gauge function $\ga:Y\to R\cup \{-\oo,\oo\}$ attached to
$(H,\Lb)$, under the model of (c03).
Further, letting $(X,d)$ be a metric space, denote (again) by
$(\succeq)$ the quasi-order on $X\times Y$ introduced
as in (d03); and finally, let $A$ be some (nonempty)
part of $X\times Y$.
The question to be posed is the same as in Section 4;
to solve it, we list the needed conditions.
The former of these is (again) (d04);
which also writes
\bit
\item[(e01)]
$P_Y(A)\incl H$\ \ (i.e.:\ $\wt y=0$ in that condition).
\eit
[For, otherwise, passing to the subset
$A_0=\{(x,y)\in X\times Y; (x,y+\wt y)\in A\}$,
this requirement is fulfilled, via
$P_Y(A_0)=P_Y(A)-\wt y$].
As a consequence, $\inf[\ga(P_Y(A)]\ge 0$
(wherefrom $\ga$ is bounded below on $P_Y(A)$).
However, the alternative $\ga[P_Y(A)]=\{\oo\}$ cannot be avoided;
so, we must accept (as a second condition)
\bit
\item[(e02)]
$P_Y(A)\cap \Dom(\ga)\ne \es$\
($\ga(y)< \oo$, for some $y\in P_Y(A)$).
\eit
A useful characterization of these is to be realized via the
composed function
$\Phi(x,y)=\ga(y), (x,y)\in X\times Y$\ (i.e.: $\Phi=\ga \circ P_Y$).
Precisely, let again $\Phi$ denote the restriction of this function
to $A$; then, (e01)+(e02) may be written as
\bit
\item[(e03)]
$\inf[\Phi(A)]\ge 0$\ and\  
$\Dom(\Phi):=\{(x,y)\in A; \Phi(x,y)< \oo\}$ is nonempty.
\eit
Now, the last condition to be imposed is a variant of (a04) above:
\bit
\item[(e04)]
each $(\succeq)$-ascending $e$-Cauchy sequence
$((x_n,y_n))\incl \Dom(\Phi)$\\
is bounded above in $A$ (modulo $(\succeq)$).
\eit
Here, $e$ stands for the semi-metric over $X\times Y$
introduced as in (d07).
\sk

We are now in position to state the second main result
in this exposition.

\btheorem \label{t3}
Let the precise conditions be in force.
Then, for each $(x_0,y_0)\in \Dom(\Phi)$ there exists
$(\bar x,\bar y)\in \Dom(\Phi)$ with the properties (\ref{101}) and
\beq \label{501}
\mbox{
$(\bar x,\bar y)\succeq (x',y')\in A$ imply 
$\bar x =x'$, $\ga(\bar y)=\ga(y')$.
}
\eeq
\etheorem

\bproof
We claim that Proposition \ref{p2} is applicable to
$(Z=X\times Y,\succeq)$, $M=A$ and $\vphi=\Phi$.
In fact, by the remarks above
$(x_1,y_1)\succeq (x_2,y_2)$ implies
$\Phi(x_1,y_1)\ge \Phi(x_2,y_2)$;
i.e., $\Phi$ is $(\succeq)$-decreasing. On the other hand,
(e03) is just (b03) (with $\vphi$ substituted by $\Phi$).
Finally, (e04) implies (b04) (with $\vphi=\Phi$); and this will
establish our claim. In fact, let  $((x_n,y_n))$ be an ascending
sequence in $\Dom(\Phi)$; i.e.,
\bit
\item[]
$\Lb(d(x_n,x_m))\le_K y_n-y_m$,\ \ whenever $n\le m$.
\eit
Combining with the subtractivity of the gauge function (cf. Section 3) yields
$$  \mbox{
$d(x_n,x_m)\le \ga(y_n-y_m)\le \ga(y_n)-\ga(y_m)$,\ \ whenever\ $n\le m$.
}
$$
The (real) sequence $(\ga(y_n))$ is descending and bounded
(by the choice of our data); hence a Cauchy one. This, added to
the above, shows that $(x_n)$ is a  $d$-Cauchy sequence in
$P_Y(A)$; or, equivalently, that $((x_n,y_n))$ is $e$-Cauchy in
$A$; wherefrom (by (e04)) the claim follows.
By Proposition \ref{p2} we therefore derive that, for
$(x_0,y_0)\in \Dom(\Phi)$ there exists
$(\bar x,\bar y)\in \Dom(\Phi)$ with the properties (\ref{101}) and
$$  \mbox{
$(\bar x,\bar y)\succeq (x',y')\in A$ imply
$\Phi(\bar x,\bar y)=\Phi(x',y')$.
}
$$
The relation in the left member of this implication yields
(see the remarks above):
$(x',y')\in \Dom(\Phi), d(\bar x,x')\le \ga(\bar y)-\ga(y')$.
Moreover, the relation in the right member of the same is just:
$\ga(\bar y)=\ga(y')$; so that (combining these)
$d(\bar x,x')=0$; wherefrom $\bar x=x'$.
This proves (\ref{501}) as well; and concludes the argument.
\eproof

As before, when $\Lb$ is taken as in (a05) (and $d$ is complete)
Theorem \ref{t3} is nothing but the related statement in
Turinici \cite{turinici-02},
obtained via similar techniques.
[Moreover, when $H$ is taken as in (d05), we have (cf. Section 3)
$$  \mbox{
$\Dom(\ga)= H$;\ hence $\Dom(\Phi)=A$ (in view of (e01));
}
$$
and Theorem \ref{t3} reduces to Theorem \ref{t2} above.
But, in the general (nonlinear) setting, this is not true].
Further aspects (of locally convex nature) may be found in
Isac and Tammer \cite{isac-tammer-04};
see also
Rozoveanu \cite{rozoveanu-00}.
For different structural extensions of these we refer to
Khanh \cite{khanh-88}.

\section{Further aspects}
\setcounter{equation}{0}

Our main results (Theorem \ref{t2} and Theorem \ref{t3})
were especially designed to 
enlarge (in two distinct manners)
the product variational principle in
Goepfert, Tammer and Z\u{a}linescu \cite{goepfert-tammer-zalinescu-00}.
Unfortunately, neither of these can extend in a direct way
the related variational statements in 
Bao and Mordukhovich
\cite[Theorem 3.4]{bao-mordukhovich-10};
for, e.g., the quasi-boundedness assumption imposed 
(by the authors) upon $x\mapsto A(x)$ is weaker than the 
boundedness condition (b04) used here.
On the other hand, the  cited results cannot extend 
Theorem \ref{t2} or Theorem \ref{t3};
because the (a02)-type  product quasi-order of the authors
is the linear version (cf. (a05)) of the "non-linear" 
product quasi-order (d03) used by us. Having these precise, 
it is natural to ask whether a common extension of all these 
variational principles is available; 
we conjecture that the answer is positive.



\begin{thebibliography}{99}

\bibitem{altman-82}
{M. Altman},
\it A generalization of the Brezis-Browder principle on ordered sets,
\rm Nonlinear Analysis, 6 (1982), 157-165.


\bibitem{bao-mordukhovich-10}
{T. Q. Bao} and {B. S. Mordukhovich}, 
\it Relative Pareto minimizers for multiobjective problems: existence and optimality conditions, 
\rm Math. Program., 122 (2010), 301-347.


\bibitem{brezis-browder-76}
{H. Brezis} and {F. E. Browder},
\it A general principle on ordered sets in nonlinear functional analysis,
\rm Advances Math., 21 (1976), 355-364.


\bibitem{carja-ursescu-93}
{O. C\^{a}rj\u{a}} and {C. Ursescu},
\it The characteristics method for a first order partial differential
equation,
\rm  An. St. Univ. "A. I. Cuza" Ia\c{s}i (S. I-a, Mat.), 39 (1993), 367--396.


\bibitem{cristescu-77}
{R. Cristescu},
\it Topological Vector Spaces,
\rm  Noordhoff Intl. Publishers, Leyden (The Netherlands), 1977.


\bibitem{ekeland-79}
{I. Ekeland},
\it Nonconvex minimization problems,
\rm Bull. Amer. Math. Soc. (New Series), 1 (1979), 443-474.


\bibitem{gerth(tammer)-weidner-90}
{C. Gerth (Tammer)} and {P. Weidner},
\it Nonconvex separation theorems and some applications in
vector optimization,
\rm J. Optim. Theory Appl., 67 (1990), 297-320.


\bibitem{goepfert-riahi-tammer-zalinescu-03}
{A. Goepfert}, {H. Riahi}, {C. Tammer} and {C. Z\u{a}linescu},
\it Variational Methods in Partially Ordered Spaces,
\rm [Canad. Math. Soc. Books in Math. vol. 17],
Springer, New York, 2003.


\bibitem{goepfert-tammer-zalinescu-00}
{A. Goepfert}, {C. Tammer} and {C. Z\u{a}linescu},
\it On the vectorial Ekeland's variational principle
and minimal points in product spaces,
\rm Nonlinear Analysis, 39 (2000), 909--922.


\bibitem{isac-96}
{G. Isac},
\it  The Ekeland's principle and Pareto $\veps$-efficiency,
\rm in "Multi--Objective Programming and Goal Programming" (M. Tamiz ed.),
pp. 148--163, L. Notes Econ. Math. Systems vol. 432, Springer, Berlin, 1996.


\bibitem{isac-tammer-04}
{G. Isac} and {C. Tammer},
\it Nuclear and full nuclear cones in product spaces:
Pareto efficiency and an Ekeland type variational principle,
\rm Positivity, 14 (2004), 1-28.


\bibitem{kang-park-90}
{B. G. Kang} and {S. Park},
\it On generalized ordering principles in nonlinear analysis,
\rm Nonlinear Analysis, 14 (1990), 159-165.


\bibitem{khanh-88}
{P. Q. Khanh},
\it On Caristi-Kirk's theorem and Ekeland's variational principle
for Pareto extrema,
\rm Bull. Polish Acad. Sci. (Math.), 37 (1988), 33-39.



\bibitem{nemeth-86}
{A. B. Nemeth},
\it A nonconvex vector minimization problem,
\rm Nonlinear Analysis, 10 (1986), 669-678.


\bibitem{rozoveanu-00}
{P. Rozoveanu},
\it Ekeland's variational principle for vector valued functions,
\rm Math. Reports [St. Cerc. Mat.], 2(52) (2000), 351-366.


\bibitem{tammer-92}
{C. Tammer},
\it A generalization of Ekeland's variational principle,
\rm Optimization, 25 (1992), 129-141.


\bibitem{turinici-89}
{M. Turinici},
\it Metric variants of the Brezis-Browder ordering principle,
\rm Demonstr. Math., 22 (1989), 213-228.


\bibitem{turinici-02}
{M. Turinici},
\it Minimal points in product spaces,
\rm An. \c{S}t. Univ. "Ovidius" Constan\c{t}a (Ser. Math.),
10 (2002), 109-122.


\bibitem{turinici-06}
{M. Turinici},
\it Brezis-Browder principles in separable ordered sets,
\rm Libertas Math., 26 (2006), 15-30.


\end{thebibliography}
\end{document}